\documentclass[12pt]{amsart}
\usepackage[pdftex]{graphicx}
\usepackage{float}
\usepackage{thumbpdf}
\usepackage{amsfonts,amssymb,amsmath,type1cm}
\usepackage[pdftex,unicode]{hyperref}

\newtheorem{theorem}{Theorem}

\newcommand{\RR}{{\mathbb R}}
\newcommand{\ZZ}{{\mathbb Z}}

\newcommand{\CC}{{\mathbb C}}

\newcommand{\pt}{{\mathcal{PT}}}

\usepackage[usenames,dvipsnames]{color}

\begin{document}

\title[On perturbations of the spectrum of $\pt$-operator]{On perturbations of the spectrum of 
one-dimensional $\pt$-symmetric periodic Schr\"odinger operator}
\author{P.G. Grinevich}
\address{Steklov Mathematical Institute of Russian Academy of Sciences\\ L.D. Landau Institute for Theoretical Physics of Russian Academy of Sciences \\ Lomonosov Moscow State University}
\email{pgg@landau.ac.ru}

\author{I.A. Taimanov}
\address{Novosibirsk State University, Novosibirsk, Russia}
\email{iskander.taimanov@mail.ru}

\date{}
\maketitle

\begin{abstract}
For  $\pt$-symmetric periodic Schr\"odinger operator, which is a small perturbation of the zero potential, we calculate the spectrum and the divisor of zeroes of the Bloch function in the leading order of the perturbation theory. In particular, we show that the analogs of lacunae of the Bloch spectrum are ellipses, and their focal points coincide with the branch points of the spectral curve. 
\end{abstract}




\vskip3mm

\hfill{To the memory of A.V. Borisov (1965--2021)}

\vskip5mm

\section{Introduction}

In the recent years a new class of differential operators attracted an attention of physicists
(see, for instance, \cite{BB,B19}).
These are so-called operators with $\pt$-symmetry.
For one-dimensional Schr\"odinger operators
\begin{equation}
\label{eq:sch1}
L = -\frac{d^2}{dx^2} + u(x), 
\end{equation}
that means that the potential $u(x)$ satisfies the condition
\begin{equation}
\label{eq:pt}
u(x) = \overline{u(-x)}.
\end{equation}

The spectral theory of periodic $\pt$-symmetric Schr\"odinger operators demonstrates many interesting features \cite{V20,V25}.

In \cite{T24} the second author (I.A.T.) described finite gap (or finite-zone) Schr\"odinger operators
with $\pt$-potentials in terms of the algebro-geometric data for the inverse spectral problem for periodic potentials. 

In the present article we address two problems:

1) how new gaps in the spectral curves are transformed under perturbations of the potential?

2) what are the analogs of the lacunae for $\pt$-potentials?

We mean by the analogs of the lacunae of the Bloch spectrum the projections of the solutions to the Dubrovin equations (see \eqref{eq:dubr} below) onto the $E$-plane where $E$ is the the eigenvalue of $L$. For real-valued potentials such defined lacunae coincide with the classical ones.

\section{Preliminary facts}

Finite gap (or finite-zone) 
\footnote{Novikov who introduced these operators in \cite{N74} coined them finite-zone because his approach was motivated by the solid state physics where the gaps in the Bloch spectrum are called ``zones of instability'' and have important physical meaning.}
potentials of the one-dimensional Schr\"odin\-ger operators with  periodic potentials 
\begin{equation}
\label{eq:sch1}
L = -\frac{d^2}{dx^2} + u(x), \ \ u(x+T) = u(x),
\end{equation}
were introduced in \cite{N74}. The inverse problem for such real-valued potentials was solved in \cite{D75,IM75}. 

We recall that for periodic potentials  on the two-dimensional space of solutions
to the equation
\[
L \psi = E\psi, \ \ E = \mathrm{const} \in \CC,
\]
there is defined the monodromy operator
\[
T(E) \psi(x) = \psi(x+T).
\]
Its eigenvalues satisfy the equation
\begin{equation}
\label{eq:blocheig}
\lambda^2 - 2r(E) \lambda +1 =0 
\end{equation}
where $r(E)$ is an entire function.
Hence the boundary points of the periodic and antiperiodic spectra $\{E\}$ are given by the equation
\[ 
r^2(E) = 1
\]
and the operator is called finite-zone if this equation has finitely many simple roots
\[
E_1,\dots,E_{2g+1}.
\]
There are also infinitely many double roots, the resonant points of the spectrum.

The eigenfunctions of $T(E)$ are called the (Floquet--) Bloch functions and the corresponding eigenvalues 
are called the Bloch multipliers. 

The (classical) spectrum of a periodic Schr\"odinger operator consists of such $E$ for which the Bloch functions are bounded, or, equivalently, $|\lambda|=1$ where $\lambda$ satisfies \eqref{eq:blocheig}.   

In \cite{D75,IM75} it is shown that for the finite gap potentials $u(x)$ the Bloch functions are parameterized 
by the points of the Riemann surface $\Gamma$ of the form
\begin{equation}
\label{eq:surface}
w^2 = R(E) = (E-E_1)\dots(E-E_{2g+1})
\end{equation}
which is called the spectral curve of the operator $L$.

Moreover, given a fixed point $x_0 \in \RR$, there is uniquely defined a global Bloch function
$\Psi(x,x_0,P)$, where $P \in \Gamma$, such that

1) it is normalized by the condition $\Psi(x_0,x_0,P) = 1$;

2) for any given $P_0 \in \Gamma$ the function $\Psi(x,x_0,P_0)$ is the Bloch function corresponding to
$(E,w) \in \Gamma$;

3) $\Psi(x,x_0,P) \approx e^{i\sqrt{E} (x-x_0)}$ as $E \to \infty$;

4) the function $\Psi$ is expressed in terms of algebraic functions corresponding to $\Gamma$.  

For real-valued potentials these facts are proved in \cite{D75,IM75} and such proofs can be generalized for complex-valued potentials. For detailed exposition of the finite-gap theory for complex-valued algebro-geometric potentials we refer to \cite{BG}. 

The analytical formulas for $u(x)$ are different:

1) in \cite{D75} the real-valued potential $u(x)$ takes the form
\[
u(x) = -2\sum_{k=1}^g \gamma_k(x) + \sum_{j=1}^{2g+1} E_j,
\]
where $\gamma_k(x) \in [E_{2k},E_{2k+1}], k=1,\dots, g$, and satisfy the Dubrovin equation
\begin{equation}
\label{eq:dubr}
\gamma_k^\prime(x) = \frac{-2i\sqrt{R(\gamma_k)}}{\prod_{j\neq k} (\gamma_k-\gamma_j)}.
\end{equation}
These quantities $\gamma_k(x)$ are the eigenvalues of the linear problem
\[
L \varphi = E \varphi, \ \ \varphi(x) = \varphi(x+T) = 0.
\]
Moreover, 
\begin{equation}
\label{eq:norm}
|\Psi(x,x_0,P)| = \frac{\prod_{k=1}^g (E - \gamma_k(x))}{\prod_{k=1}^g (E-\gamma_(x_0))};
\end{equation}

2) in \cite{IM75} the potential is given by the Its--Matveev formula
\begin{equation}
\label{eq:im}
u(x) = -2\frac{d^2}{dx^2} \theta(Ux + Z) +C, 
\end{equation}
where $\theta$ is the theta function of the Riemann surface \eqref{eq:surface}, $U$ 
is a constant vector and $C$ is a constant.

Of course these formulas are related via the Abel transform. 

For complex-valued potentials 

1) the formula \eqref{eq:im} is straightforwardly applied;

2) the Dubrovin equations are also valid but  the points $\gamma_k$ does not lie above the lacunae $[E_{2k},E_{2k+1}]$ of the spectrum but evolve along complicated cycles on the surface \eqref{eq:surface}.  

In \cite{T24} the second author (I.A.T.) used \eqref{eq:im} for characterizing the so-called $\pt$-potentials
which are finite-zone.

However to apply the Dubrovin approach to $\pt$-potentials it needs to know the cycles drawn by 
$\gamma_k(x)$ on the surface \eqref{eq:surface} and their description even in the $\pt$-situation is unknown.

\section{Main results}

We consider two problems:

1) how the resonant points of the spectrum are deformed under small perturbations of the potential?

2) how look like the dynamics of the points $\gamma_k(x)$ under the Dubrovin equations \eqref{eq:dubr}?

We  consider both of them in the model case when the zero potential is perturbed.
Such a method was applied to the linear problem associated with the nonlinear Schr\"odinger equation  in \cite{GS18,GS22}.
We focus our attention on the  first-order perturbation theory,
because it is very illustrative. The first-order formulas are explicit, very simple,
and they allow to demonstrate qualitatively what happens if the perturbation has
moderate strength.

These perturbations in general does not preserve the finite gap property however locally, near the resonant point, they give an answer to the first question (Theorem 1). 

Also in the first order of approximation they allow to demonstrate what are the trajectories of $\gamma_k(x)$ 
under \eqref{eq:dubr}. For real-valued potentials such trajectories are the preimages of the lacunae 
with respect to the projection 
\[
\Gamma \to \CC, \ \ (E,w) \to E.
\]
For generic $\pt$-potentials they look very differently (Theorem 2).

Without loss of generality we assume that the potential is $2\pi$-periodic:
\[
u(x+2\pi) = u(x).
\]
The $\pt$-constraint imply that 
\begin{equation}
\label{eq:sch3}
u(x)= \sum\limits_{l\in\ZZ}c_l e^{ilx}, \ \ c_l \in \RR.
\end{equation}
Without loss of generality we may assume $c_0=0$.

Let us fix a Bloch multiplier
$$
\varkappa = e^{2\pi i\alpha}.
$$
Then we have the standard Fourier basis in the space $\mathcal{L}(\varkappa)$
\begin{equation}
\label{eq:basis1}  
\psi_k = e^{i(k+\alpha)x}, \ \ \ k \in \ZZ.
\end{equation}
The restriction of $L$ to the space  $\mathcal{L}(\varkappa)$ is described by the following matrix:
\begin{equation}
\label{eq:sch4}
\left[L_{\varkappa}\right]_{mn}= \delta_{mn}(m+\alpha)^2 + c_{m-n}
\end{equation}
We consider the linearized direct spectral transform. Assume that
\begin{equation}
\label{eq:sch5}
c_j = O(\varepsilon), \ \ \varepsilon\ll 1.
\end{equation}
Let
\[
  L^{(0)} = -\frac{d^2}{dx^2}.
\]
Outside the $\varepsilon$-neighbourhoods of the resonant points 
\begin{equation}
\label{eq:sch6}
  E_n^{(0)} = \frac{n^2}{4}, \ \ n=1,2,3,\ldots.
\end{equation}
to every energy level $E$ there corresponds two different Bloch multipliers and this picture is saved by small perturbations of $L^{(0)}$. Therefore
we are left to study the spectrum of $L$ near the resonant point.

Let us calculate the perturbed spectrum near the resonant point $E_n$. Let
\begin{equation}
\label{eq:res2}
  k_n= \frac{n}{2}, \ \ \varkappa = e^{i(k_n+\delta)x} = e^{i\left( \frac{n}{2} + \delta\right) x},
\end{equation}
the near-resonant pair of eigenfunctions is
\begin{equation}
\label{eq:res3}
 \psi_{+} = e^{i(k_n+\delta)x}, \ \ \psi_{-} = e^{i(-k_n+\delta)x}.
\end{equation}

The corresponding block in the matrix representing $L$ has the form
\begin{equation}
  \label{eq:res4}
P_n=  \begin{bmatrix}
E_n^{(0)} + \delta^2 + n \delta & c_n \\ c_{-n} & E_n^{(0)} + \delta^2 - n \delta
 \end{bmatrix}   
\end{equation}
and the eigenvalues of $P_n$ are:
\begin{equation}
  \label{eq:res5}
\lambda_n^{(\pm)}= E_n^{(0)} + \delta^2 \pm \sqrt{n^2\delta^2+c_nc_{-n}}     
\end{equation}

\begin{figure}
  \centering{
    1.) $c_nc_{-n}>0$: \includegraphics[width=7cm]{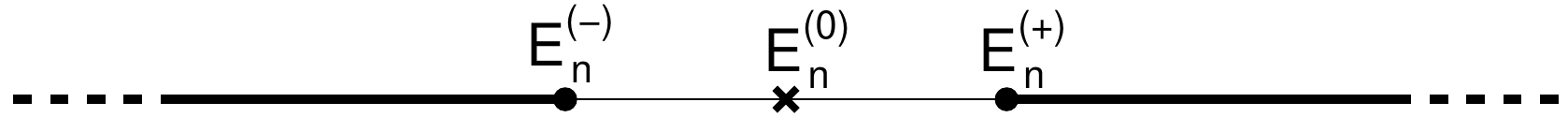}\\
    \medskip
    2.) $c_nc_{-n}=0$:\includegraphics[width=7cm]{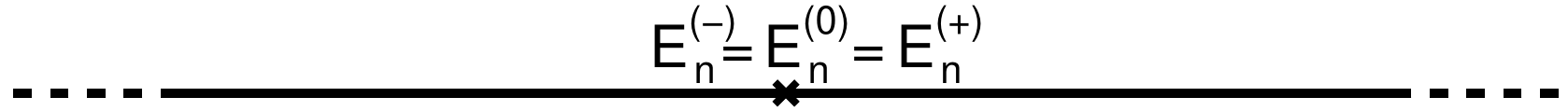}\\
    \medskip
    3.) $c_nc_{-n}<0$: \includegraphics[width=7cm]{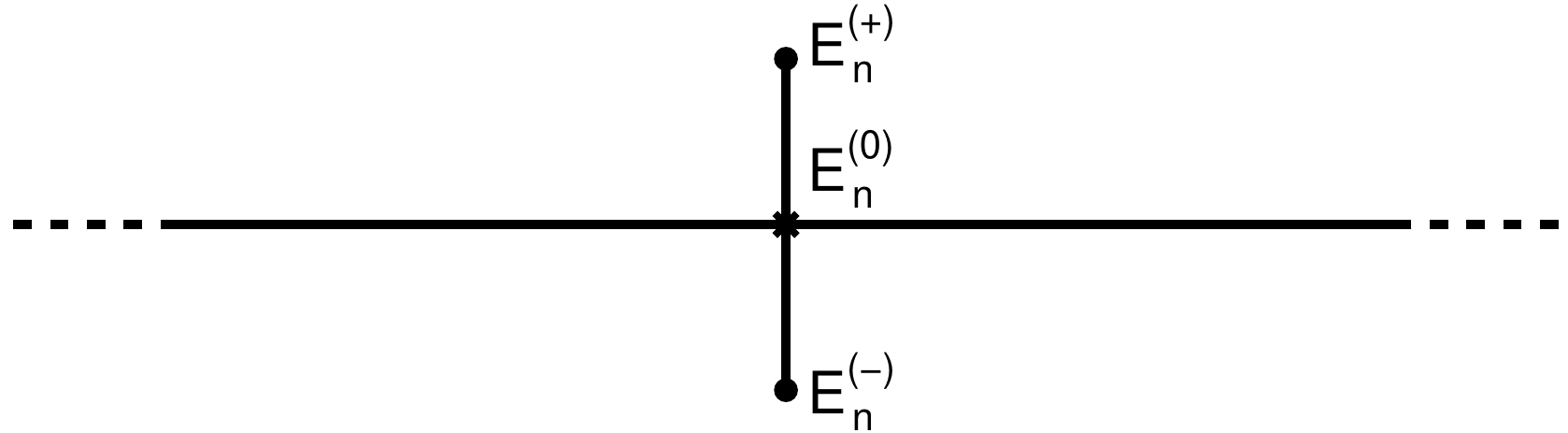}\\ 
    \caption{At the top: the case $c_nc_{-n}>0$. The perturbation generates a gap in the spectrum. In the middle:  the case $c_nc_{-n}=0$. In the leading order of the perturbation theory the spectral curve has a double point. At the bottom: the case $c_nc_{-n}<0$. The is no gap in the spectrum, moreover, the spectrum contains an interval perpendicular to the real line.\label{fig:1}} }    
\end{figure}

We have three possible configuration (see Figure~\ref{fig:1}):
\begin{enumerate}

\item $c_nc_{-n}>0$. 
In this case the perturbation generates a gap in the spectrum, the branch points are: 
\[
E_n^{(\pm)}= E_n^{(0)} \pm \sqrt{c_n c_{-n}};
\]

\item $c_nc_{-n}=0$. In this case we have no splitting of the resonant point in the leading order calculation, 
and in this order the resonant point becomes a double point on the spectral curve. 
However the existence of a multiple point in the spectrum is unstable to arbitrary small perturbations and we can make conclusion of the existence of a double point in the spectrum only after counting all orders of approximation; 

\item $c_nc_{-n}<0$. In this case the perturbation generates an additional band in the spectrum, the branch points are: 
\[
E_n^{(\pm)}= E_n^{(0)} \pm i \sqrt{|c_nc_{-n}|}
\]
and the approximated spectrum is parameterized by real valued parameter $\delta$ such that
\[ 
|\delta| \leq \frac{|c_n c_{-n}|}{n^2}
\]
for which 
\[
|\lambda^{\pm}_n| \leq 1.
\]
\end{enumerate}

\begin{theorem}
For the small $\pt$-perturbations \eqref{eq:sch3} of the zero potential
the resonant point \eqref{eq:sch6} of the spectrum of $L$ 

1) is opened into a gap for $c_n c_{-n} >0$;

2) is opened into a branch of the spectrum which is a curve transversal to the spectral half-line $\{E >0\}$ 
and  bounded by the complex conjugate (branch) points $E^{+}_n = \overline{E}^{-}_n$ for $c_n c_{-n} < 0$.  
\end{theorem}

Now we may answer the original question about the analogs of lacunae for $\pt$-potentials. 

For $c_n c_{-n} <0$ we have a branch of the spectrum which is a curve bounded by two complex-conjugate points. To such a branch there corresponds the dynamics of $\gamma(x)$ which is, by \eqref{eq:norm},
is the zero of the $\Psi(x,x_0,\gamma(x))$.

In the leading order the Bloch eigenfunction
$$
L\Psi = \lambda \Psi
$$
with the Bloch multiplier 
\[
\varkappa= e^{i\left(\frac{n}{2}+\delta\right)x},
\]
where
\[
\lambda=E_n+\delta^2+\tilde\lambda, \ \ \ 
\tilde\lambda^2=n^2\delta^2 + c_n c_{-n}
\]
can be written up to multiples as either
\begin{equation}
\label{eq:sch7}
\Psi = \bigg[-c_n\exp\left(i\frac{n}{2}x\right) + (n\delta-\tilde\lambda) \exp\left(-i\frac{n}{2}x\right)\bigg]
\exp(i\delta x)  
\end{equation}
or 
\begin{equation}
\label{eq:sch8}
\Psi = \bigg[(-n\delta-\tilde\lambda)\exp\left(i\frac{n}{2}x\right) - c_{-n} \exp\left(-i\frac{n}{2}x\right)\bigg]
\exp(i\delta x)  
\end{equation}
Therefore the condition $\Psi=0$ is equivalent to:

\begin{align*}
\label{eq:sch10}
  n\delta-\tilde\lambda &= c_n e^{i n x} \\
  -n\delta-\tilde\lambda &= c_{-n} e^{-i n x},
\end{align*}
and, finally
\begin{align*}
\tilde\lambda &= -\frac{c_n e^{i n x}+ c_{-n} e^{-i n x}}{2}\\
\delta  &= \frac{c_n e^{i n x} - c_{-n} e^{-i n x}}{2}.
\end{align*}

We conclude

\begin{theorem}
For the small $\pt$-perturbations \eqref{eq:sch3} of the zero potential with $c_n c_{-n}< 0$
in the principal order of approximation the projections $\gamma(x)$ of the zeroes of the Bloch function $\Psi$
are ellipses, and their focal point coincide with the branch points. For $x=0$ the projections of the divisor  
$\gamma(0)$ lie on the real line.
\end{theorem}

\vskip3mm

The work was supported by RSCF (grant \textnumero 	24-11-00281).

\end{document}